\documentclass{amsart}
\usepackage{amssymb, latexsym}
\usepackage{verbatim, enumerate}
\DeclareMathOperator{\im}{Im}

\DeclareMathOperator{\re}{Re}

\newcommand{\eps}{{\varepsilon}}
\newcommand{\D}{{\mathbb D}}

\newcommand{\C}{{\mathbb C}}
\newcommand{\R}{{\mathbb R}}

\newcommand{\E}{{\mathbb E}}

\newcommand{\Szego}{{Szeg\H{o} }}
\newcommand{\su}{{\mathbb U}}
\newcommand{\so}{{\mathbb S \mathbb O}}
\newcommand{\eql}{{\phantom{E}\stackrel{\text{law}}{=}\phantom{E}}}

\newcommand{\interval}{{[-2,2]}}

\theoremstyle{plain}
\newtheorem{thm}{Theorem}
\newtheorem{lemma}[thm]{Lemma}
\newtheorem{remark}[thm]{Remark}
\newtheorem{prop}[thm]{Proposition}

\numberwithin{equation}{section} \numberwithin{thm}{section}

\begin{document}

\title[Linear statistics of point processes]{Linear statistics of point processes via orthogonal polynomials}
\author{E. ~Ryckman}
\address{253-37 Mathematics, California Institute of Technology, Pasadena CA 91125, USA}
\email{eryckman@caltech.edu}

\begin{abstract}
For arbitrary $\beta > 0$, we use the orthogonal polynomials techniques developed in \cite{kncmv, Killip-Nenciu:matrix-models} to study certain linear statistics associated with the circular and Jacobi $\beta$ ensembles.  We identify the distribution of these statistics then prove a joint central limit theorem.  In the circular case, similar statements have been proved using different methods by a number of authors.  In the Jacobi case these results are new.
\end{abstract}

\maketitle

%%%%%%%%%%%%%%%%%%%%%%%%%%%%%%%%%%%%%%%%%%%%%%%%%%%
%
%
%                                   Section
%
%
%%%%%%%%%%%%%%%%%%%%%%%%%%%%%%%%%%%%%%%%%%%%%%%%%%%

\section{Introduction and main results}
In this paper we will study two families of random point processes, the circular and Jacobi $\beta$ ensembles.  For $\beta > 0$, the circular $\beta$ ensemble with $n$ points, $C\beta E(n)$, is the random point process on the unit circle where  for any symmetric function $f$ we have
\begin{equation}\label{beta measure}
\E_n^\beta (f) = \frac{1}{Z_n^\beta} \int_{-\pi}^\pi \dots \int_{-\pi}^\pi f(e^{i\eta_1}, \dots, e^{i\eta_n})   |\Delta(e^{i\eta_1},\dots,e^{i \eta_n})|^\beta \frac{d\eta_1}{2\pi}\dots\frac{d\eta_n}{2\pi}.
\end{equation}
Here, $\Delta$ is the usual Vandermonde determinant
$$\Delta (z_1, \dots , z_n) = \prod_{1 \leq i < j \leq n} (z_j - z_i)$$ and the partition function is $$Z_n^\beta = \frac{\Gamma(\frac{n \beta}{2} + 1)}{[\Gamma(\frac{\beta}{2} + 1)]^n}$$ as was shown in \cite{MR0258644, Killip-Nenciu:matrix-models, MR0144627}.  The particular cases $\beta = 1, 2$, and $4$ correspond to the classical orthogonal, unitary, and symplectic circular ensembles (COE, CUE, and CSE).  In particular, if $\su(n)$ is the set of $n \times n$ unitary matrices and $U \in \su (n)$ is chosen uniformly with respect to Haar measure, then its eigenvalues are a sample from $C (\beta = 2) E$.  It is also known that the above measure is the Gibbs measure for an $n$-partical Coulomb gas at inverse temperature $\beta$, but we will not pursue this viewpoint here.  For more background see, for instance, \cite{MR2129906}.

The second point process we consider is the Jacobi $\beta$ ensemble of $n$ points, $J\beta E_{a,b}(n)$.  Here, the points $x_1, \dots ,x_n$ are confined to live on $\interval$ with joint probability density proportional to
\begin{equation}\label{jacobi}
|\Delta(x_1, \dots, x_n)|^\beta \prod_j (2-x_j)^a (2+x_j)^b dx_1 \dots dx_n
\end{equation}
where $\beta > 0$ and $a, b > -1$.  The specific case $\beta = 2, a=b=-1/2$ corresponds to (half) the eigenvalues of  a matrix chosen randomly from $\so (2n)$, see \cite{Killip-Nenciu:matrix-models}.  In particular, any statement about $J\beta E$ yields a corresponding statement about $\so (2n)$, and we will not formulate these results separately.

Given a sample of points $\{ e^{i\eta_1}, \dots , e^{i\eta_n} \}$ from $C\beta E$, we will consider the statistic $\prod (1 - e^{-i(\theta - \eta_k)}).$   Note that if these points are the eigenvalues of some unitary matrix $U$, then the above statistic is just the characteristic polynomial $\det (1 - e^{-i\theta}U)$.  In the Jacobi case we will instead consider $\prod (E - x_k)$.  Our first result is to identify the distribution of these statistics.  For this, recall that $X$ is Beta-distributed on $[0,1]$ with parameters $s$ and $t$, written $X \sim B_{[0,1]}(s,t)$, if it has the density function $$\frac{\Gamma(s+t)}{\Gamma(s)\Gamma(t)} x^{s-1}(1-x)^{t-1}dx.$$  Similarly we will write $X \sim B_{[-1,1]}(s,t)$ if the corresponding density function is $$\frac{2^{1-s-t}\Gamma(s+t)}{\Gamma(s)\Gamma(t)}(1-x)^{s-1}(1+x)^{t-1}dx.$$  Note that in this case $\E(X) = \frac{t-s}{t+s}$ and $\E(X^2) = \frac{(t-s)^2 + (t+s)}{(t+s)(t+s+1)}$.

\begin{thm}\label{main 1}
We have the following statements concerning distributions:
\newline

\noindent $(C\beta E)$:  Let $\beta > 0$ be arbitrary and say $e^{i\eta_1}, \dots, e^{i\eta_n}$ be a sample from $C\beta E $.  Then for all $\theta \in \R$ $$\prod_{k=1}^n (1 - e^{-i(\theta-\eta_k)}) \eql \prod_{k=1}^n (1 + e^{i\theta_k}\sqrt{X_{k-1}})$$ where $\theta_i, X_j$ are all independent, each $\theta_k$ is uniformly distributed on $S^1$ and $X_{j} \sim B_{[0,1]}\bigl(1, \frac{j \beta}{2}\bigr)$.\newline

\noindent $(J\beta E)$:  Let $\beta > 0$ be arbitrary and say $x_1, \dots, x_n \in \interval$ are a sample from $J\beta E_{a,b}$.  Then $$\prod_{k=1}^n (\pm 2 - x_k) \eql 2 \prod_{k=0}^{2n-2} (1 - (\pm 1)^{k+1 }X_{k})$$ where the $X_k$'s are all independent and distributed as
$$X_k \sim \begin{cases}B_{[-1,1]}\Biggl(\frac{k}{4}\beta + a + 1, \frac{k}{4}\beta + b + 1 \Biggr), \quad \text{k even}\\ B_{[-1,1]}\Biggl(\frac{k-1}{4}\beta + a + b + 2, \frac{k + 1}{4}\beta\Biggr), \quad \text{k odd.}\end{cases}$$
\end{thm}

\begin{remark}
In the $C\beta E$ case, rotation invariance shows that the distribution of $\prod (1 - e^{-i(\theta - \eta_k)})$ is independent of $\theta$.  On the real line no such invariance is present, and we have been unable to find a simple expression of $\prod (E - x_k)$ for $E \neq \pm 2$.  (However, see Proposition 6.1 of \cite{Killip-Nenciu:matrix-models}, where the expectation for general $E$ is expressed in terms of Jacobi polynomials.)
\end{remark}

Our second main result concerns the limiting behavior of our statistics as $n$ tends to infinity.

\begin{thm}\label{jl}
The following asymptotic joint laws hold:\newline

\noindent $(C \beta E)$: Let $\beta > 0$ be arbitrary, $e^{i\eta_1}, \dots, e^{i\eta_n}$ be a sample from $C\beta E$, and $Z_n (\theta) = \sum_{k=1}^n \log  (1 - e^{-i(\theta-\eta_k)})$.  Then for any distinct $e^{i\theta_1},  \dots, e^{i\theta_M}$ in $S^1$, as $n \rightarrow \infty$ the joint law of
\begin{equation*}
\frac{1}{\sqrt{\log n}} \Bigl(  Z_n(\theta_1), \dots,  Z_n(\theta_M) \Bigr)
\end{equation*}
converges to the joint law of $M$ independent complex normal random variables with mean zero and variance $1/\beta$.\newline

\noindent $(J \beta E)$:  Let $\beta > 0$ be arbitrary, $x_1, \dots, x_n$ be a sample from $J\beta E_{a,b}$, and $Z_n (E) = \sum_{k=1}^n \log(E - x_k)$.  Then for any distinct $e^{i\theta_1}, \dots , e^{i\theta_M}$ in $S^1$, as $n \rightarrow \infty$ the joint law of 
\begin{equation*}
\frac{1}{\sqrt{\log n}} \Bigl( Z_n (2\cos \theta_1) - E_n(\theta_1), \dots, Z_n(2 \cos \theta_M)- E_n(\theta_M)    \Bigr)
\end{equation*}
converges to the joint law of $M$ independent complex normal random variables with mean zero and variance $1/\beta$, where $$E_n(\theta) = (C_0\delta_0(\theta) + C_\pi \delta_\pi(\theta))\log n - i n \theta$$
\begin{equation}\label{C}
C_0 = \frac{2a+1}{\beta}-\frac{1}{2} \quad\quad\quad C_\pi = \frac{2b+1}{\beta} - \frac{1}{2}
\end{equation}
and $\delta_\theta$ is a Dirac measure at $\theta$.
\end{thm}

\begin{remark}
While the appearance of delta measures at the edge of the interval may strike some readers as odd, it is not unprecedented in this context.  See, for instance, formula (3.54) in \cite{MR1487983}.
\end{remark}

Other probabilistic statements are certainly possible.  For instance, a straightforward adaptation of the arguments in \cite{bhny} yields iterated logarithm laws for the various ensembles considered above.  We do not pursue this point further.

In the circular case, results of this type are not new and go back at least to the papers \cite{MR1794267, MR1794265} of Keating and Snaith.  Motivated by connections to analytic number theory, they studied the distribution of $Z_n (\theta)$ when $U$ is chosen from $C \beta E$, $SO(2n)$, or $USp(2n)$.  Using the Mellin-Fourier transform they identified the distribution of $Z_n$ by calculating the moment generating function averaged over the relevant collection of matrices.  From this one may readily deduce a version of Theorem \ref{jl} (but only in the circular case, and only in the one-point case $M=1$).

Motivated in part by this work, Bourgade et al. \!\!\! \cite{bhny, bnr} used probabilistic methods to study $\su (n)$ and $\so (2n)$.  They proved the $C\beta E$ case of Theorem \ref{main 1}, and from this characterization they were able to recover the (one-point) central limit theorem of Keating and Snaith, as well as develop some new results concerning the speed of convergence (in particular the iterated logarithm law mentioned earlier).  

The question of joint asymptotic behavior was first addressed by Hughes, Keating, and O'Connell in \cite{MR1844632}, where they used Szeg\H{o}-type limit theorems to deduce asymptotic properties of the joint law.  In particular they proved the circular case of Theorem \ref{jl} and also deduced a variety of large-deviation results.  (We should note here that while the above papers only considered $C \beta E$ for the particular values $\beta = 1, 2, 4$, their calculations easily extend to handle all $\beta > 0$.)

All the previously-existing results discussed so far pertain only to the circular case.  There are relatively fewer results concerning point processes on the real line.  We mention \cite{MR2092032} wherein the moment generating function for the real part of the log of the characteristic polynomial is expressed for many ensembles (including the circular and Jacobi ensembles) in terms of Jack polynomials and generalized hypergeometric functions.  As far as we can tell though, the orthogonal polynomials technique--and resulting joint central limit theorem for the Jacobi case that we establish--do not appear anywhere in the current literature.

The remainder of this paper is organized as follows: Section \ref{op section} collects some basic results we will need concerning orthogonal polynomials.  We prove Theorem \ref{main 1} in Section \ref{distribution section} and Theorem \ref{jl} in Section \ref{jl section}.  

It is a pleasure to thank Rowan Killip for many helpful conversations during the preparation of this work, as well as the referee for relevant references to the literature.

%%%%%%%%%%%%%%%%%%%%%%%%%%%%%%%%%%%%%%%%%%%%%%%%%%%
%
%
%                                   Section
%
%
%%%%%%%%%%%%%%%%%%%%%%%%%%%%%%%%%%%%%%%%%%%%%%%%%%%

\section{Orthogonal polynomials background}\label{op section}
Throughout this paper we will make use of some results from the theory of orthogonal polynomials on the unit circle and real line, or put another way, the study of CMV and Jacobi matrices.  In this section we collect the results we will need without proof (for a more detailed discussion see \cite{Simon:OPUC1,Simon:OPUC2} and \cite{Teschl:book}).  

Let $d\mu$ be a probability measure on $S^1$ that consists of $n$ point-masses.  The monomials $1, z, \dots, z^{n-1}$ form a basis for $L^2(d\mu)$, so by applying the Gram-Schmidt procedure we arrive at an orthogonal basis of monic polynomials $\Phi_k(z) = z^k + \dots$, $k = 0, \dots, n-1$.  Define the reversed polynomials by 
\begin{equation}\label{reversed}
\Phi_k^\ast (z) = z^k \overline{\Phi_k (1/\bar{z})}.
\end{equation}
The first important fact is that these polynomials obey the \Szego recurrence equation:
\begin{align}\label{szego recurrence}
\begin{split}
\Phi_{k+1}(z) &= z \Phi_k(z) - \overline{\alpha_k} \Phi_k^\ast (z)\\
\Phi_{k+1}^\ast (z) &= \Phi_k^\ast (z) - z \alpha_k \Phi_k (z)
\end{split}
\end{align}
for some sequence $\alpha_k$, termed the Verblunsky coefficients associated to $d\mu$.  It is easy to see that for $k < n-1,$ $\alpha_k \in \D$ ($=$ the unit disk in $\C$) and $\alpha_{n-1} \in S^1$.   

Given a sequence of Verblunsky coefficients let $\rho_k = \sqrt{1 - |\alpha_k|^2}$ and define $$\Xi_k = \begin{bmatrix}\overline{\alpha_k} & \rho_k \\ \rho_k & - \alpha_k \end{bmatrix}$$ for $0 \leq k \leq n-2$, and $\Xi_{-1} = [1]$, $\Xi_{n-1} = [\overline{\alpha_{n-1}}]$.  From these form the $n \times n$ block-diagonal matrices $$L_n = \text{diag}(\Xi_0, \Xi_2, \Xi_4, \dots)\quad\quad\quad M_n=\text{diag}(\Xi_{-1}, \Xi_1, \Xi_3, \dots)$$ and then the CMV matrices 
\begin{equation}\label{cmv}
C_n = L_nM_n.
\end{equation}

\begin{remark}\label{associated}
We now have four pieces of data: the measure $d\mu$, the orthogonal polynomials $\Phi_k$, the Verblunsky coefficients $\alpha_k$, and the CMV matrix $C_n$.  It is known that any one of these sets uniquely determines the others.  So, for instance, we may speak of the orthogonal polynomials associated to a set of Verblunsky coefficients, and we will do so throughout this paper.
\end{remark}

The eigenvalues of CMV matrices will provide the link between orthogonal polynomials and point processes.  The first step is:

\begin{prop}\label{det poly relation}
$\Phi_n(z) = \det(z - C_n)$ and $\Phi_n^\ast (z) = \det (1 - z \overline {C_n})$ where $(\overline C)_{ij} = \overline{C_{ij}}$.  In particular, $\det (1 - C_n) = \Phi_n(1) = \overline{\Phi_n^\ast (1)}$.
\end{prop}

We now add some randomness to the CMV construction.  If $\nu > 1$ we will say that a complex-valued random variable Z with values in $\D$ is $\Theta_\nu$-distributed if its density function is $$\frac{\nu - 1}{2} (1 - |z|^2)^{(\nu - 3)/2}\frac{d^2 z}{\pi}.$$  We will say $Z \sim \Theta_1$ if $Z$ is uniformly-distributed on $S^1$.  With this notation we have:

\begin{thm}[\cite{Killip-Nenciu:matrix-models}]\label{kn}
Let $\beta > 0$ and let $\alpha_k \sim \Theta_{\beta k+1}$ for $0 \leq k \leq n-1$ be independent.  Then the CMV matrices $C_n$ give a matrix model for $C \beta E$.  That is, their eigenvalues are distributed according to the circular ensemble \eqref{beta measure}.
\end{thm}

We now move to the real line as follows: given a measure $d\mu$ on the unit circle with $d\mu(\overline z) = d\mu(z)$, define a measure $d\nu$ on $\interval$ by $$\int_{-\pi}^\pi f(2\cos \theta) d\mu(\theta) = \int_{-2}^2 f(x) d\nu(x).$$  As before, Gram-Schmidt yields a sequence $P_k(x) = x^k + \dots$ of monic orthogonal polynomials, and then normalized polynomials $p_k = \frac{P_k}{\|P_k\|}$.  These obey a 3-term recurrence equation $$x p_k(x) = a_{k+1}p_{k+1}(x) + b_{k+1}p_k(x) + a_k p_{k-1}(x)$$ for some sequences $a_k > 0$, $b_k \in \R$.  From these sequences we build the (cut-off) Jacobi matrix $$J_n = \begin{bmatrix} b_1 & a_1 & & \\ a_1 & b_2 & \ddots & &\\ &\ddots & \ddots & a_{n-1}\\ & & a_{n-1}  & b_n\end{bmatrix}.$$  Again we have four pieces of data (measure, polynomials, recurrence coefficients, matrix), and each uniquely determines the others.

It is a famous observation of \Szego \cite{Szego:book} that the orthogonal polynomials for $d\nu$ and $d\mu$ are related
\begin{equation}\label{line circle}
P_k(z + z^{-1}) = \frac{z^{-k}\Phi_{2k}(z) + z^k \Phi_{2k}(z^{-1})}{1 - \alpha_{2k-1}}
\end{equation}
and as before we have a relation between $J_n$ and the orthogonal polynomials:

\begin{prop}\label{det poly relation 2}
$\det (x - J_n) = P_n(x)$ and $\det (\pm 2 - J_n) = \frac{2}{1 - \alpha_{2n-1}}\Phi_{2n}(\pm 1)$.
\end{prop}

Again we introduce some randomness:

\begin{thm}[\cite{Killip-Nenciu:matrix-models}]\label{knj}
Let $\alpha_k$, $0 \leq k \leq 2(n-1)$ be independent and distribued according to
\begin{equation}\label{knj dist}
\alpha_k \sim \begin{cases}B_{[-1,1]}\Biggl(\frac{k}{4}\beta + a + 1, \frac{k}{4}\beta + b + 1 \Biggr), \quad \text{k even}\\ B_{[-1,1]}\Biggl(\frac{k-1}{4}\beta + a + b + 2, \frac{k + 1}{4}\beta\Biggr), \quad \text{k odd.}\end{cases}
\end{equation}
Let $\alpha_{2n-1} = \alpha_{-1} = -1$, and define
\begin{align}\label{geronimus}
\begin{split}
a_{k+1}^2 &= (1 - \alpha_{2k-1})(1 - \alpha_{2k}^2)(1 + \alpha_{2k+1})\\
b_{k+1} &= (1 - \alpha_{2k-1}) \alpha_{2k} - (1 + \alpha_{2k-1}) \alpha_{2k-2}
\end{split}
\end{align}
for $0 \leq k \leq n-1$.  Then the Jacobi matrices $J_n$ give a matrix model for $J \beta E_{a,b}$.  That is, the eigenvalues of $J_n$ are distributed according to the Jacobi ensemble \eqref{jacobi}.
\end{thm}

One final remark: strictly speaking, Theorems \ref{kn} and \ref{knj} differ from those appearing in \cite{Killip-Nenciu:matrix-models} by a relabeling of the $\alpha$'s.  That this does not change the above statements is the content of Proposition B.2 in \cite{Killip-Nenciu:matrix-models}, and the rotation-invariance of the $\Theta$-distribution.  We chose this presentation merely to simplify our formulas.

%%%%%%%%%%%%%%%%%%%%%%%%%%%%%%%%%%%%%%%%%%%%%%%%%%%
%
%
%                                   Section
%
%
%%%%%%%%%%%%%%%%%%%%%%%%%%%%%%%%%%%%%%%%%%%%%%%%%%%

\section{Distribution Results}\label{distribution section}
In this section we prove Theorem \ref{main 1}.  Let us begin with $C \beta E$.  By the rotation-invariance of \eqref{beta measure}, the distribution of $\prod_{k=1}^n (1 - e^{-i(\theta-\eta_k})$ is independent of $\theta$, so without loss we let $\theta=0$.  In this case $\prod_{k=1}^n (1 - e^{i\eta_k}) = \Phi_n(1)$ by Proposition \ref{det poly relation}.  We will now relate $\Phi_n(1)$ to the Verblunsky coefficients, and then the Verblunsky coefficients to the random variables appearing in the theorem.

\begin{lemma}\label{1}
Let the $\alpha$'s be distributed as in Theorem \ref{kn}, and let $\Phi_k$ be the associated orthogonal polynomials.  Then $$\Phi_n(1)  \eql \prod_{k=0}^{n-1} (1 + \alpha_k).$$
\end{lemma}

\begin{proof}
From the definition of the reversed polynomials, $\Phi_k^\ast (1) = \overline{\Phi_k(1)}$ for all $k$.  Using this and \Szego recursion we get $$\Phi_n (1) = \Phi_{n-1}(1) - \overline{\alpha_{n-1}} \Phi_{n-1}^\ast(1) = \Phi_{n-1}(1) - \overline{\alpha_{n-1}} \overline{\Phi_{n-1}(1)}.$$
As the distribution of the $\alpha$'s is rotationally symmetric, this is equal in law to $\Phi_{n-1}(1) (1 + \alpha_{n-1})$.  The result now follows by iteration.
\end{proof}

\begin{lemma}\label{2}
Let the $\alpha$'s be distributed as in Theorem \ref{kn}.  Then $$\alpha_k \eql e^{i\theta} \sqrt{X_{k}}$$ where $\theta$ and $X_k$ are all independent, $\theta$ is uniformly distributed on $[0,2\pi)$ and $X_{k} \sim B_{[0,1]}\bigl(1, \frac{\beta k}{2}\bigr)$.
\end{lemma}

\begin{proof}
$X_{k}$ has as its density function $$\frac{\beta k}{2}(1-x)^{\frac{\beta k}{2}-1}\chi_{[0,1]}(x) dx.$$

By Theorem \ref{kn}, $\alpha_k$ has as its density function
\begin{equation*}
\frac{\beta k+1-1}{2}(1-|z|^2)^{\frac{\beta k+1-3}{2}} \chi_\D (z) \frac{d^2 z}{\pi} = 
\frac{\beta k}{2}(1-|z|^2)^{\frac{\beta k}{2}-1}\chi_\D (z)\frac{d^2 z}{\pi}.
\end{equation*}
This is rotationally-symmetric, hence the $e^{i\theta}$ term.  The square root comes from $|z|^2$ vs $x$.
\end{proof}

\begin{proof}[Proof of Theorem \ref{main 1} $(C \beta E)$]
By Theorem \ref{kn} it suffices to consider only eigenvalues of CMV matrices.  In this case the result follows from  Proposition \ref{det poly relation} and Lemmas \ref{1} and \ref{2}.
\end{proof}

We now turn to the $J \beta E$ part of Theorem \ref{main 1}.  Essentially the same proof as that of Lemma \ref{1} shows:

\begin{lemma}\label{3}
If the $\alpha$'s are distributed as in Theorem \ref{knj}, then $$\Phi_{2n}(\pm 1) \phantom E = \phantom E 2 \prod_{k=0}^{2n-2} (1 - (\pm 1)^{k+1} \alpha_k).$$
\end{lemma}

\begin{proof}[Proof of Theorem \ref{main 1} $(J \beta E_{a,b})$]
Combine Proposition \ref{det poly relation 2}, Theorem \ref{knj}, and Lemma \ref {3}.
\end{proof}

%%%%%%%%%%%%%%%%%%%%%%%%%%%%%%%%%%%%%%%%%%%%%%%%%%%
%
%
%                                   Section
%
%
%%%%%%%%%%%%%%%%%%%%%%%%%%%%%%%%%%%%%%%%%%%%%%%%%%%

\section{Joint laws}\label{jl section}
In this section we prove Theorem \ref{jl}.  In both the circular and Jacobi cases the result follows from some basic estimates and a version of the central limit theorem for martingales (see Proposition \ref{mclt} below).  We begin by setting up some notation.  Let  $z=e^{i\theta}$ and $$B_k(z) = z \frac{\Phi_k (z)}{\Phi_k^\ast (z)} = e^{i \psi_k (\theta)}$$ where the random continuous functions $\psi_k$ are defined recursively by $$\psi_0 (\theta) = 0 \quad\quad \psi_{k+1}(\theta) = \psi_k (\theta) + \theta -2 \im \Upsilon (\psi_k(\theta), \alpha_k)$$ $$\Upsilon(\psi_k(\theta),\alpha_k) =  \log(1 - \alpha_k e^{i\psi_k (\theta)}).$$  In particular notice that $$\psi_k(0) = 0 \quad\quad\text{and}\quad\quad \psi_k(\pi) = \begin{cases} 0 \text{ , $k$ odd}\\ \pi \text{ , $k$ even}. \end{cases}$$

Next define
\begin{gather*}
\widetilde \Upsilon (\psi, \alpha) = \begin{cases} -[\alpha - \E(\alpha)]e^{i\psi} , \text{ Jacobi}\\ -\alpha e^{i\psi}, \text{\phantom{iiiiiiiiiiiiiii}Circular} \end{cases}\\
T(n,\theta) = \sum_{k=0}^{n-1} \re \Upsilon (\psi_k (\theta), \alpha_k) \quad\quad \widetilde T(n,\theta) = \sum_{k=0}^{n-1} \re \widetilde \Upsilon (\psi_k (\theta), \alpha_k)\\
S(n,\theta) = \sum_{k=0}^{n-1} \im \Upsilon (\psi_k (\theta), \alpha_k) \quad\quad \widetilde S(n,\theta) = \sum_{k=0}^{n-1} \im \widetilde \Upsilon (\psi_k (\theta), \alpha_k).
\end{gather*}

\begin{lemma}\label{st relation}
Let the $\alpha$'s be distributed as in Theorem \ref{kn} (the circular case).  Then $$\re Z_n (\theta) = T(n,\theta) \quad\quad \im Z_n (\theta) = S(n,\theta).$$  If instead the $\alpha$'s are distributed as in Theorem \ref{knj} (the Jacobi case), then $$\re Z_n(2 \cos \theta) = T(2n,\theta) \quad\quad \im Z_n(2 \cos \theta) = S(2n,\theta) - in\theta.$$
\end{lemma}

\begin{proof}
We first consider the $C \beta E$ case.  From \Szego recursion \eqref{szego recurrence} we find 
\begin{equation}\label{phistar recur}
\Phi_{n+1}^\ast (z) = \Phi_n^\ast (z) \bigl( 1 - \alpha_n B_n (z) \bigr) = \prod_{k=0}^n \bigl( 1 - \alpha_k B_k (z) \bigr).
\end{equation}
So by Proposition \ref{det poly relation}, $\im Z_n (\theta) =  \im \log \Phi_n^\ast (\theta) =  S(n,\theta)$ and similarly $\re Z_n (\theta) = T(n, \theta)$.  

Now consider the $J \beta E$ case.  From the relation \eqref{line circle} we have 
$$(1 - \alpha_{2n-1}) P_n (2 \cos \theta) = e^{-in\theta} \Phi_{2n}^\ast (\theta) \Biggl( 1 + \frac{\Phi_{2n}(\theta)}{\Phi_{2n}^\ast (\theta)} \Biggr).$$  
From \Szego recursion \eqref{szego recurrence} and the fact that $\alpha_{2n-1} = -1$ we get $\Phi_{2n}(\theta) = \Phi_{2n}^\ast (\theta)$, so by \eqref{phistar recur} the above equation becomes
$$P_n(2 \cos \theta) = e^{-i n \theta} \Phi_{2n}^\ast (\theta) = e^{-i n \theta} \prod_{k=0}^{2n-1} (1 - \alpha_k e^{i\psi_k(\theta)}).$$  So by Proposition \ref{det poly relation 2} we have $$Z_n (2 \cos \theta) = \log  P_n (2 \cos \theta) = -in\theta + \sum_{k=0}^{2n-1} \log (1 - \alpha_k e^{i \psi_k(\theta)})$$ as claimed.
\end{proof}

We now turn to the asymptotic behavior of $S$ and $T$, where our analysis will follow that of \cite{killip:gaussian} quite closely.  By the recursive definition of $\psi_k$ and that $\E (\widetilde \Upsilon) = 0$, we see that $\widetilde S(n,\theta)$ and $\widetilde T(n,\theta)$ are martingales with respect to the sigma algebras $$\mathcal M_k = \sigma (\alpha_0, \dots, \alpha_{k-1}).$$  This will allow us to use the following version of the central limit theorem (one should think of $\Psi$ as essentially the real or imaginary part of $\Upsilon$):

\begin{prop}\label{mclt}
Fix $\theta_1, \dots, \theta_M \in [0, 2\pi)$ distinct and let $\widetilde \Psi(\psi_k(\theta), \alpha_k)$ be a martingale with respect to $\mathcal M_k$.  Suppose that there is $\Psi(\psi_k (\theta), \alpha_k)$ so that as $n \rightarrow \infty$ 
\begin{equation}\label{cond1}
\frac{1}{\log n} \sum_{k=0}^{n-1} \E \bigl( \widetilde \Psi (\psi_k (\theta_j), \alpha_k) \widetilde \Psi (\psi_k(\theta_l),\alpha_k) | \mathcal M_k \bigr) \stackrel{L^1}{\rightarrow} \sigma^2 \delta_{jl}
\end{equation}
\begin{equation}\label{cond2}
\frac{1}{\log^2 n} \sum_{k=0}^{n-1} \E \bigl( | \widetilde \Psi (\psi_k (\theta_j),\alpha_k)|^4 \bigr) \rightarrow 0
\end{equation}
\begin{equation}\label{cond3}
\frac{1}{\sqrt{\log n}} \E \Biggl( \Biggl| \sum_{k=0}^{n-1} \Psi(\psi_k(\theta_j), \alpha_k) - \widetilde \Psi (\psi_k(\theta_l), \alpha_k) \Biggr| \Biggr) \rightarrow 0
\end{equation}
for all $1 \leq j,l \leq M$.  Then the random variables $$R(n,\theta_j) = \frac{\sum_{k=0}^{n-1}\Psi(\psi_k(\theta_j),\alpha_k)}{\sqrt{\log n}}$$ converge to independent normal random variables with mean zero and variance $\sigma^2$.
\end{prop}

This is proved in \cite{killip:gaussian}.  Basically, the first two conditions allow one to apply the usual Martingale Central Limit Theorem (\cite{MR1852999} Chapter 6) to $$\widetilde R (n, \theta_j) = \sum_{k=0}^{n-1}\widetilde \Psi (\psi_k(\theta_j),\alpha_k)$$ while the last condition relates convergence of $R$ to that of $\widetilde R$.

In the $C \beta E$ case we'll apply this result to the real and imaginary parts of $\Upsilon$ and $\widetilde \Upsilon$, then make an additional argument to show that these are asymptotically independent.  The $J \beta E$ case is essentially the same, but slightly more involved at the edges of the interval $\interval$. The next two lemmas summarize the estimates we'll need to verify conditions \eqref{cond1}--\eqref{cond3}.

\begin{lemma}[$C\beta E$]
Suppose $\phi,\psi \in \R$ and $\alpha \sim \Theta_\nu$.  Let $\widetilde \chi$ be either $\re \widetilde \Upsilon$ or $\im \widetilde \Upsilon$ and let $\chi$ be the same but without the tilde.  Then
\begin{equation}\label{a}
\E \bigl(\widetilde \chi (\psi, \alpha) \widetilde \chi (\phi, \alpha) \bigr) = \frac{2 \cos (\psi-\phi)}{\nu + 1}
\end{equation}
\begin{equation}\label{b}
\E \bigl( |\widetilde \chi (\psi,\alpha)|^4 \bigr) \lesssim \frac{1}{\nu^2}
\end{equation}
\begin{equation}\label{c}
\E \bigl( | \chi (\psi, \alpha) - \widetilde \chi (\psi, \alpha) | \bigr) \lesssim \frac{1}{\nu^2}.
\end{equation}
\begin{equation}\label{d}
\E \bigl(\re \widetilde \Upsilon (\psi, \alpha) \im \widetilde \Upsilon (\phi, \alpha) \bigr) = \frac{\sin (\phi-\psi)}{2(\nu + 1)}
\end{equation}
where all the implicit constants are independent of $\alpha, \nu, \psi, \phi$.
\end{lemma}

\begin{proof}
Consider \eqref{d}.  Using the definition of $\Theta_\nu$ we have 
\begin{align*}
\E \bigl(\re \widetilde \Upsilon (\psi, \alpha) \im \widetilde \Upsilon (\phi, \alpha) \bigr) &= \frac{\nu -1}{2} \int\!\!\!\int_\D \re (z e^{i\psi}) \im (z e^{i\phi}) (1 - |z|^2)^{\frac{\nu-3}{2}}\frac{d^2 z}{\pi}\\
&= \frac{\nu -1}{2} \int_0^{2\pi} \int_0^1 r^3 \cos(\theta+\psi) \sin(\theta + \phi) (1 - r^2)^{\frac{\nu-3}{2}} dr \frac{d\theta}{2\pi}\\
&= \frac{1}{2(\nu+1)}\sin (\phi - \psi).
\end{align*}

All the other proofs are similar and we do not present them.  (For the reader seeking more details, see \cite{killip:gaussian} where the case $\widetilde \chi = \im \widetilde \Upsilon$, $\chi  = \im \Upsilon$ appears as Lemma 2.5.  The proofs for the real parts are exactly the same.)
\end{proof}

\begin{lemma}[$J \beta E$]
Suppose $\phi,\psi \in \R$ and $\alpha \sim B_{[-1,1]}(s,t)$.  Then
\begin{equation}\label{e}
\E \bigl( \re \widetilde \Upsilon(\psi,\alpha) \re  \widetilde \Upsilon (\phi, \alpha) \bigr) = \frac{2st}{(s+t)^2(s+t+1)} \bigl( \cos (\psi-\phi) + \cos (\psi + \phi) \bigr)
\end{equation}
\begin{equation}\label{e2}
\E \bigl( \im \widetilde \Upsilon(\psi,\alpha) \im  \widetilde \Upsilon (\phi, \alpha) \bigr) = \frac{2st}{(s+t)^2(s+t+1)} \bigl( \cos (\psi-\phi) - \cos (\psi + \phi) \bigr)
\end{equation}
\begin{equation}\label{e3}
\E \bigl( \re \widetilde \Upsilon(\psi,\alpha) \im  \widetilde \Upsilon (\phi, \alpha) \bigr) = \frac{2st}{(s+t)^2(s+t+1)} \bigl( \sin (\psi+\phi) - \sin (\psi - \phi) \bigr).
\end{equation}
If the $\alpha$'s are distributed as in Theorem \ref{knj}, and $\widetilde \chi$ is either of $\re \widetilde \Upsilon$ or $\im \widetilde \Upsilon$, then
\begin{equation}\label{f}
\E \bigl(| \widetilde \chi (\psi,\alpha_k)|^4 \bigr) \lesssim \frac{1}{(k+1)^2}
\end{equation}
\begin{equation}\label{h}
\E \bigl(| \widetilde \chi (\psi,\alpha_k)| \bigr) \lesssim \frac{1}{(k+1)^{1/2}}.
\end{equation}
Finally if $\theta \neq 0, \pi$ then
\begin{equation}\label{g}
\E \bigl( |\re \Upsilon (\psi_k (\theta), \alpha_k) + \alpha_k \cos (\psi_k(\theta)) + \frac{1}{2}\alpha_k^2 \cos(2 \psi_k(\theta))  | \bigr) \lesssim \frac{1}{(k+1)^2}
\end{equation}
while for all $\psi$
\begin{equation}\label{h}
\E \bigl( |\im \Upsilon (\psi, \alpha_k) + \alpha_k \sin (\psi) + \frac{1}{2}\alpha_k^2 \sin(2 \psi)  | \bigr) \lesssim \frac{1}{(k+1)^2}.
\end{equation}
\end{lemma}

\begin{proof}
The proof is a straightforward adaptation of that just presented for the circular case, so we do not present the details (but see Lemma 2.6 in \cite{killip:gaussian}, which presents the proof for the imaginary parts).
\end{proof}

We will need one more result:
\begin{lemma}\label{sumbyparts}
Given real-valued sequences $\eps_k, X_k, Y_k$ with $X_{k+1} = X_k + \delta + Y_k$ for some $\delta \in (0,2\pi)$,
$$
\Biggl| \sum \eps_k e^{iX_k} \Biggr| \leq \frac{2 \|\eps_k\|_{\ell^{\infty}} + \| \eps_k - \eps_{k-1}\|_{\ell^1} + \| \eps_k Y_k \|_{\ell^1}}{| 1 - e^{i\delta}|}.
$$
\end{lemma}

\begin{proof}
This appears as Lemma 2.7 in \cite{killip:gaussian} and is proved there (it is essentially just summation by parts).
\end{proof}

\begin{proof}[Proof of Theorem \ref{jl}]
We begin with $C \beta E$.  The first step is to show that \eqref{cond1}--\eqref{cond3} are satisfied for $\widetilde \Psi = \re \widetilde \Upsilon$ or $\widetilde \Psi = \im \widetilde \Upsilon$.  Since the proofs are identical we only present the case $\re \widetilde \Upsilon$ (one can find the proof for the imaginary parts in \cite{killip:gaussian}).

\underline{\eqref{cond1}}:  Using \eqref{a} can rewrite the left-hand side of \eqref{cond1} as 
$$\frac{1}{\log n} \sum_{k=0}^{n-1} \frac{1}{\beta k +2} \cos(\psi_k(\theta_j) - \psi_k(\theta_l))$$ which immediately settles the case $j=l$ with $\sigma^2 = \tfrac{1}{\beta}$.  The case $j \neq l$ follows from Lemma \ref{sumbyparts} with 
\begin{gather*}
X_k = \psi_k(\theta_j) - \psi_k(\theta_l) \quad\quad Y_k = -2 \im \Upsilon (\psi_k(\theta_j),\alpha_k) +2  \im \Upsilon (\psi_k(\theta_l), \alpha_k)\\
\eps_k = \frac{1}{\beta k +2}\quad\quad\delta = \theta_j - \theta_l.
\end{gather*}

\underline{\eqref{cond2}}: This follows directly from \eqref{b}.

\underline{\eqref{cond3}}: Using $\E (\re \Upsilon) = \E (\re \widetilde \Upsilon) = 0$ and \eqref{c} we can write
\begin{align*}
\frac{1}{\log n} \E \Biggl( \Biggl| \sum_{k=0}^{n-1} \re \Upsilon(\psi_k(\theta),&\alpha_k) - \re \widetilde \Upsilon (\psi_k(\theta),\alpha_k) \Biggr|^2 \Biggr)\\ 
&= \frac{1}{\log n} \sum_{k=0}^{n-1} \E \Bigl( | \re \Upsilon (\psi_k(\theta),\alpha_k) - \re \widetilde \Upsilon (\psi_k(\theta),\alpha_k) |^2 \Bigr)\\
& \lesssim \frac{1}{\log n} \sum_{k=0}^{n-1}\frac{1}{(k+1)^2} \rightarrow 0.
\end{align*}
Condition \eqref{cond3} then follows from Cauchy-Schwarz.

Now that \eqref{cond1}--\eqref{cond3} are verified, Proposition \ref{mclt} shows that the random variables $T(n,\theta_1), \dots, T(n, \theta_j)$ converge to independent normal random variables with mean zero and variance $1/\beta$.  A similar statement holds with $S$ replacing $T$.  To show that the limits of the $T$'s are independent of the limits of the $S$'s (and thus $T + i S$ converge to independent complex normals), we must show that the potential contributions from the cross-terms $\re \widetilde \Upsilon(\psi_k(\theta_j),\alpha_k) \im \widetilde \Upsilon (\psi_k(\theta_l), \alpha_k)$ vanish asymptotically for all $j, l$.  That is, we must show
\begin{equation}\label{cross}
\frac{1}{\log n}\sum_{k=0}^{n-1} \E \bigl( \re \widetilde \Upsilon(\psi_k(\theta_j),\alpha_k) \im \widetilde \Upsilon (\psi_k(\theta_l), \alpha_k) | \mathcal M_k \bigr) \stackrel{L^1}{\rightarrow} 0.
\end{equation}

By \eqref{d} we see the left-hand side above is (up to an irrelevant constant)
$$\frac{1}{\log n} \sum_{k=0}^{n-1} \frac{\sin \bigl( \psi_k(\theta_l) - \psi_k(\theta_j) \bigr)}{k+1} .$$  This already settles the case $j=l$.  The case $j \neq l$ follows from Lemma \ref{sumbyparts} with 
\begin{gather*}
X_k = \psi_k(\theta_j) - \psi_k(\theta_l) \quad\quad Y_k = -2 \im \Upsilon (\psi_k(\theta_j),\alpha_k) +2  \im \Upsilon (\psi_k(\theta_l), \alpha_k)\\
\eps_k = \frac{1}{ k+1}\quad\quad\delta = \theta_j - \theta_l.
\end{gather*}
By the comments in the previous paragraph this finishes the proof for $C\beta E$.

We now turn to $J\beta E$.  By Lemma \ref{st relation} we have $$Z_n (2 \cos \theta) - E_n(\theta) = \bigl[ T(2n,\theta) - \bigl( C_0 \delta_0 (\theta) + C_\pi \delta_\pi (\theta) \bigr) \log n \bigr] + i S(2n,\theta).$$  Again we must verify the conditions \eqref{cond1}--\eqref{cond3} hold, and show asymptotic independence of the real and imaginary parts.  The arguments for \eqref{cond1} and \eqref{cond2} are a straightforward adaptation of the proof just presented, but now using \eqref{e}--\eqref{h}.  Similarly, the proof of asymptotic independence follows that given in the $C\beta E$ case and uses \eqref{e3}.  We omit these easy calculations. This leaves us to check condition \eqref{cond3}, which we will only prove for the real part (the proof of the imaginary part is easier, and can also be found in \cite{killip:gaussian}).

First suppose that $\theta \neq 0, \pi$.  In this case we use \eqref{g} to rewrite the left-hand side of \eqref{cond3} as $$\frac{1}{\sqrt{\log n}} \E \Biggl( \Biggl| \sum_{k=0}^{2n-1} \E (\alpha_k) \cos (\psi_k (\theta)) - \frac{1}{2} \alpha_k^2 \cos (2 \psi_k (\theta)) \Biggr| \Biggr) + O \Bigl( \frac{1}{\sqrt{\log n}} \Bigr).$$  From Lemma \ref{sumbyparts} with $\delta = \theta, X_k = \psi_k, \eps_k = \E (\alpha_k)$ we find 
$$\E\Biggl( \Biggl| \sum_{k=0}^{2n-1} \E(\alpha_k) \cos(\psi_k(\theta)) \Biggr| \Biggr) = O(1)$$ and similarly $$\quad  \E \Biggl( \Biggl| \sum_{k=0}^{2n-1} \E (\alpha_k^2) \cos (2 \psi_k (\theta)) \Biggr| \Biggr) = O(1).$$ The desired result then follows from combining these two estimates with Cauchy-Schwarz and the bound
\begin{align*}
\E \Biggl( \Biggl| \sum_{k=0}^{2n-1} \bigl( \alpha_k^2 - \E(\alpha_k^2) \bigr) \cos(2\psi_k(\theta))\Biggr|^2 \Biggr) &= \sum_{k=0}^{2n-1} \E \Bigl( \bigl( \alpha_k^2 - \E (\alpha_k^2)\bigr)^2 \cos^2 (2 \psi_k (\theta)) \bigr) \Bigr)\\
&\leq \sum_{k=0}^{2n-1} \E (\alpha_k^4) \lesssim \sum_{k=0}^{2n-1} \frac{1}{(k+1)^2}
\end{align*}
(which follows easily from \eqref{knj dist} and basic properties of the Beta distribution).

We now turn to the cases $\theta = 0, \pi$.  Here will need the following asymptotics (which again are easily verified using \eqref{knj dist} and basic properties of the Beta distribution): 
\begin{gather}\label{sum asymptotics}
\begin{split}
\Biggl( \Biggl|\sum_{k=0}^{n-1}\E (\alpha_{2k+1}) - \frac{\beta - 2(a+b+2)}{2\beta}\log n \Biggr|  
+ \Biggl|\sum_{k=0}^{n-1}\E (\alpha_{2k}) - \frac{b-a}{\beta}\log n \Biggr| +\\
\Biggl|\sum_{k=0}^{n-1}\E(\alpha_{2k+1}^2) - \frac{1}{\beta}\log n \Biggl|+ \Biggl|\sum_{k=0}^{n-1}\E(\alpha_{2k}^2) - \frac{1}{\beta}\log n \Biggr| \Biggr)\log^2 n = o(1).
\end{split}
\end{gather}

Now consider
\begin{align*}
\E \Biggl( \Biggl| & \sum_{k=0}^{2n-1}  \re \Upsilon(\psi_k(\!\!\phantom{b}^0_\pi), \alpha_k) - \re \widetilde \Upsilon (\psi_k (\!\!\phantom{b}^0_\pi ), \alpha_k) -  C_{\!\!\phantom{b}^0_\pi}  \log n    \Biggr| \Biggr) =\\ &\frac{1}{2}\E \Biggl( \Biggl| \sum_{k=0}^{n-1} (\alpha_{2k+1}^2+\alpha_{2k}^2) + 2 \bigl(\E (\alpha_{2k+1} \pm \E (\alpha_{2k}) \bigr) + 2 C_{\!\!\phantom{b}^0_\pi} \log n \Biggr| \Biggr) + O(1).
\end{align*}
We must verify that this is $o(\sqrt{\log n})$.  If we expand
\begin{align*}
\E \Biggl( \Biggl| \sum_{k=0}^{n-1} &\bigl( \alpha_{2k+1}^2 + \alpha_{2k}^2 \bigr)  + 2\bigl(\E (\alpha_{2k+1} \pm \E (\alpha_{2k}) \bigr) + 2C_{\!\!\phantom{b}^0_\pi} \log n \Biggr|^2 \Biggr).
\end{align*}
in powers of $C$ and use \eqref{sum asymptotics}, we find it has the form $$\Bigl( 4C^2 + 4C R + R^2 \Bigr) \log^2n + o(\log n)$$ where $R = \frac{\beta-2[2+(b+a)\pm(b-a)]}{\beta}$.  With the choice of $C_{\!\!\phantom{b}^0_\pi}$ given in \eqref{C} this is $o(\log n)$, and the final result follows by Cauchy-Schwarz.
\end{proof}

%%%%%%%%%%%%%%%%%%%%%%%%%%%%%%%%%%%%%%%%%%%%%%%%%%%
%
%
%                                 Bibliography
%
%
%%%%%%%%%%%%%%%%%%%%%%%%%%%%%%%%%%%%%%%%%%%%%%%%%%%

\bibliographystyle{plain}
\bibliography{bibliography}

\end{document}